\documentclass[titlepage,12pt]{article}
\usepackage{amssymb}
\usepackage{amsfonts}
\begin{document}

{\Large \bf Junction Conditions, \\ \\ Resolution of Singularities and \\ \\
Nonlinear Equations of Physics } \\ \\

{\bf Elem\'{e}r E Rosinger} \\
Department of Mathematics \\
and Applied Mathematics \\
University of Pretoria \\
Pretoria \\
0002 South Africa \\
eerosinger@hotmail.com \\ \\

{\bf Abstract} \\

For large classes of systems of polynomial nonlinear PDEs necessary and sufficient conditions
are given for the existence of solutions which are discontinuous across hyper-surfaces. These
PDEs contain the Navier-Stokes equations, as well as those of General Relativity and
Magneto-Hydrodynamics. \\ \\

{\bf 1. Preliminary Remarks} \\

There has for longer been an interest in finding {\it junction condition} across
hyper-surfaces of discontinuities for solutions of various nonlinear equations of Physics,
such as in General Relativity, Magneto-Hydrodynamics, and so on. Two approaches in this regard
have been pursued in the literature. One of them is trying to keep the equations and introduce
hard to deal with integral formulations of the respective conditions, while the other is
introducing ad-hoc simplifying assumptions which are dictated by mathematical convenience,
rather than physical reasons. \\

In view of the obvious disadvantages of both such approaches, in Braunss, a direct method was
suggested for a certain limited class of nonlinear equations, aiming to obtain the desired
junction conditions from the formulation of {\it weak solution conditions} across the
respective discontinuities. \\
However, in Braunss, this approach faces the following problem which is not addressed in a
rigorous way in view of the fact that the equations dealt with are {\it nonlinear}. Namely,
due to the discontinuities across the hyper-surfaces, the weak solutions include Heaviside
functions. Thus in the respective equations one becomes involved with nonlinear operations,
among them products with Heaviside functions and their derivatives, the Dirac distributions,
as well as possibly with the derivatives of Dirac distributions. Consequently, the respective
weak solution conditions fail to have a rigorous formulation in terms of the customary and
essentially linear theories of Schwartz distributions or Sobolev spaces. \\
Furthermore, the class of nonlinear PDEs dealt with in Braunss is rather limited. \\

This direct method of Braunss in formulating the junction conditions, a method with an
obviously sound and worthwhile underlying idea, was approached rigorously for the first time
in Rosinger [6, pp. 139-162], see also Rosinger [7,8], where it was applied to the large class
of arbitrary {\it polynomial nonlinear} PDEs, which thus contain among others the
Navier-Stokes equations, as well as those of General Relativity and Magneto-Hydrodynamics. \\
Such a rigorous approach was made possible by the use of the differential algebras of
generalized functions first introduced and developed in Rosinger [1-15], see also Mallios \&
Rosinger [1-3], Rosinger \& Walus [1,2], Mallios [1,2]. \\

One of the consequences is the identification of a very large class of nonlinear systems of
PDEs, called {\it resoluble}, for which such junctions conditions exist. \\
Furthermore, an {\it explicit characterization} for the existence of such junction conditions
is obtained for that class of nonlinear systems of PDEs. \\

{\it Resolution of singularities} in the above context means, therefore, the rigorous
expression of the conditions required upon solutions in the mentioned algebras across
hyper-surfaces of discontinuities. \\

The main results in this regard are reviewed here, and the respective proofs can be found in
Rosinger [6-8]. \\
Certain improvements and clarifications in the formulation of the earlier results are
presented. \\

Related to the above, one should note the following less well known, yet nontrivial fact.
Certain possible restrictive conditions imposed on the classes of nonlinear systems of PDEs
may happen to come from Physics. Indeed, as shown by Rubel for ODEs, and by Buck for PDEs, and
mentioned briefly in section 10, there are such equations which may seem to have in a certain
sense too many solutions. \\ \\

{\bf 2. The Class of Polynomial Nonlinear Systems of PDEs} \\

The polynomial nonlinear systems of PDEs dealt with are of the form \\

(2.1)~~~ $ \sum_{\,1 \,\leq\, i \,\leq\, h_\beta}~ \left \{~ c_{\beta,\, i} ( x )~ \left [~
               \prod_{\, 1 \,\leq\, j \,\leq\, k_{\beta,\, i}}~
                  D^{p_{\beta,\, i,\, j}}~U_{\alpha_{\beta,\, i,\, j}}~ ( x )
                   ~\right ] ~\right \}
                        ~=~ f_\beta ( x ) $ \\

where $x \in \Omega \subseteq \mathbb{R}^n$, with $\Omega$ nonvoid, open. The unknown
functions are $U = ( U_1,~\ldots~, U_a ) : \Omega \longrightarrow \mathbb{R}\,^a$, while $1
\leq \alpha_{\beta,\, i,\, j} \leq a,~ 1 \leq \beta \leq b,~ p_{\beta,\, i,\, j} \in
\mathbb{N}^n$, and the functions $c_{\beta,\, i},~ f_\beta \in {\cal C}^\infty ( \Omega )$ are
given. \\

The system of PDEs (2.1) is called of {\it type (MH)}, if and only if each of the associated
partial differential operators \\

(2.2)~~~ $ T_\beta ( x, D ) U ( x ) ~=~
                 \sum_{\,1 \,\leq\, i \,\leq\, h_\beta} \left \{ c_{\beta,\, i} ( x ) \left [
                   \prod_{\, 1 \,\leq\, j \,\leq\, k_{\beta,\, i}}
                     D^{p_{\beta,\, i,\, j}}~U_{\alpha_{\beta,\, i,\, j}} ( x )
                                 \right ] \right \} $ \\

can be written in the form \\

(2.3)~~~ $ \begin{array}{l}
                T_\beta ( x, D ) U ( x ) ~=~ \\ \\
                 ~~~=~ \sum_{\,1 \,\leq\, \rho \,\leq\, r_\beta}~
                         L_{\, \beta,\, \rho} ( x, D )~
                           \left [~ U_{\alpha_{\, \beta,\, \rho}}~ ( x )
                             P_{\, \beta,\, \rho} ( x, D )~
                               U_{\alpha_{\, \beta,\, \rho}}~ ( x ) ~\right ]
            \end{array} $ \\

where $L_{\, \beta,\, \rho} ( x, D )$ are linear partial differential operators with ${\cal
C}^\infty$-smooth coefficients and of order $m_{\, \beta,\, \rho}$, while $P_{\, \beta,\, \rho}
( x, D )$ are linear partial differential operators with ${\cal C}^\infty$-smooth coefficients
and of order {\it at most one}. \\

It is easy to see that the Navier-Stokes equations, as well as the equations of General
Relativity and those of Magneto-Hydrodynamics are of such type (MH). \\

The condition (2.3) can be written in the equivalent form where \\

$~~~~~~ U_{\alpha_{\, \beta,\, \rho}}~ ( x ) P_{\, \beta,\, \rho} ( x, D )~ U_{\alpha_{\,
                                                                  \beta,\, \rho}}~ ( x ) $ \\

are replaced respectively with \\

(2.4)~~~ $ \sum_{\,1 \,\leq\, \alpha,\,\, \alpha\,' \,\leq\,\, a}~
                                    U_\alpha ( x )~
                             P_{\, \beta,\, \rho,\, \alpha,\, \alpha\,'}~ ( x, D )~
                               U_{\alpha\,'}~ ( x ) $ \\

here $P_{\, \beta,\, \rho,\, \alpha,\, \alpha\,'}~ ( x, D )$ being again linear partial
differential operators with ${\cal C}^\infty$-smooth coefficients and of order at most
one. \\ \\

{\bf 3. The Singularities} \\

We turn now to the class of singularities dealt with. As is well known, for every closed
subset $\Gamma \subseteq \Omega$, there exists a ${\cal C}^\infty$-smooth function $\gamma :
\Omega \longrightarrow \mathbb{R}$, such that \\

(3.1)~~~ $ \Gamma ~=~ \{~ x \in \Omega ~|~ \gamma ( x ) ~=~ 0 ~\} $ \\

Here however, we shall impose the restriction that \\

(3.2)~~~ $ \Gamma ~~\mbox{has zero Lebesgue measure} $ \\ \\

{\bf 4. The Type of Discontinuous Solutions} \\

The solutions $U = ( U_1,~\ldots~, U_a ) : \Omega \longrightarrow \mathbb{R}\,^a$ of (MH) type
systems of nonlinear PDEs (2.1) we are interested in are supposed to satisfy \\

(4.1)~~~ $ U ~~\mbox{are}~~{\cal C}^\infty-\mbox{smooth on}~~ \Omega \setminus \Gamma $ \\

thus they can only have discontinuities across the singularity sets $\Gamma$. Consequently,
we shall be looking for such solutions which are of the form \\

(4.2)~~~ $ U ( x ) ~=~
               U_- ( x ) + \left [~ U_+ ( x ) - U_- ( x ) ~\right ] H_\gamma ( x ),~~~
                                               x \in \Omega $ \\

where $U_-$ and $U_+$ are classical solutions of the PDEs (2.1) on the whole of $\Omega$,
while the Heaviside type functions $H_\gamma : \Omega \longrightarrow \mathbb{R}$ are defined
by \\

(4.3)~~~ $ H_\gamma ( x ) ~=~ \begin{array}{|l}
                                   ~~~ 0 ~~~\mbox{if}~~ \gamma ( x ) ~\leq~ 0 \\ \\
                                   ~~~ 1 ~~~\mbox{if}~~ \gamma ( x ) ~>~ 0
                               \end{array} $ \\ \\

{\bf 5. Formulation of the Problem} \\

We are looking for the {\it necessary and sufficient junction conditions} satisfied by

\begin{itemize}

\item the (MH) type nonlinear systems of PDEs (2.1),

\item the singularities $\Gamma$ and their representations $\gamma$ in (3.1),

\item the classical solutions $U_-$ and $U_+$ of the PDEs (2.1) on the whole of $\Omega$,

\end{itemize}

so that $U$ in (4.2) are solutions of the PDEs (2.1) on the whole of $\Omega$. \\ \\

{\bf 6. Characterization of the Existence of Junction Conditions} \\

In suitable differential algebras of generalized functions, one obtains \\

{\bf Theorem 6.1.} \\

The function $U$ in (4.2) is a global solution on $\Omega$ of the (MH) type nonlinear system
of PDEs (2.1), if an only if the following {\it junction conditions} are satisfied for each
$1 \leq \beta \leq b$, in the neighbourhood of the singularity set $\Gamma$, namely \\

$ \begin{array}{l}
              \sum_\rho L_{\, \beta,\, \rho} \left \{
              \left [ \sum_{\,\alpha,\,\alpha\,'} \left [
                                    ( U_+ )_\alpha
                             P_{\, \beta,\, \rho,\, \alpha,\, \alpha\,'}
                               ( U_+ )_{\alpha\,'} -
                                    ( U_- )_\alpha
                             P_{\, \beta,\, \rho,\, \alpha,\, \alpha\,'}
                               ( U_- )_{\alpha\,'} \right ] \right ] H_\gamma \right \} + \\ \\
              + ~\frac{1}{2} \sum_\rho L_{\, \beta,\, \rho} \left \{ \left [
                     \sum_{\,\alpha,\,\alpha\,'}
                         \left [ ( U_+ )_\alpha + ( U_- )_\alpha \right ]
                         \left [ ( U_+ )_{\alpha\,'} - ( U_- )_{\alpha\,'} \right ] \right ]
                             Q_{\, \beta,\, \rho,\, \alpha,\, \alpha\,'} H_\gamma
                         \right \} = \\ \\
               = f_\beta
           \end{array} $ \\

where $Q_{\, \beta,\, \rho,\, \alpha,\, \alpha\,'}$ is the first order homogeneous part of
$P_{\, \beta,\, \rho,\, \alpha,\, \alpha\,'}$. \\

{\bf Remark 6.1.} \\

In view of the way the Heaviside type functions $H_\gamma$ appear in the junction conditions
above, it is obvious that these junction conditions are trivially satisfied outside of the
singularity sets $\Gamma$, that is, on the open and dense subsets $\Omega \setminus \Gamma$. \\
Therefore, these junction conditions act only in a neighbourhood of the singularity sets
$\Gamma$. \\
The fact that their actions are not restricted to the singularity sets $\Gamma$ alone, but
also involve their neigbourhood, results from the structure of the differential algebras of
generalized functions within which these junction conditions hold. \\
It should be noted that such an involvement of a neighbourhood of singularities, and the
impossibility to reduce the issues involved only to the singularity sets is typical for
various customary theories of generalized functions, and in particular, for the Schwartz
distributions or the Sobolev spaces. For instance, in the case of the one dimensional Dirac
$\delta$ distribution, this distribution vanishes on $\mathbb{R} \setminus \{ 0 \}$. However,
its behaviour at $0 \in \mathbb{R}$ can only be described if a neighbourhood of that point is
considered, and not only that point alone. \\ \\

{\bf 7. Resoluble Systems of Nonlinear PDEs} \\

The above results can further be extended to a yet more large class of nonlinear systems of
PDEs given by \\

{\bf Definition 7.1.} \\

The polynomial nonlinear system of PDEs (2.1) is called {\it resoluble}, if and only if the
respective nonlinear partial differential operators (2.2) are of the form \\

(7.1)~~~ $ \begin{array}{l}
               T_\beta ( x, D ) U ( x ) ~=~ \\ \\
                 ~~~=~ \sum_{\,1 \,\leq\, \rho \,\leq\, r_\beta}~
                         \left [\, T_{\, \beta,\, \rho} ( x, D )~
                          \left (\, \psi ( x ), \chi ( x ) \,\right ) \,\right ]
                             D^{p_{\, \beta,\, \rho}}
                                  \left [~ \omega ( x ) ~\right ]\,^{l_{\, \beta,\, \rho}}
            \end{array} $ \\

for $x \in \Omega$ and $1 \leq \beta \leq b$, whenever \\

(7.2)~~~ $ U ( x ) ~=~ \psi ( x ) + \chi ( x )\, \omega ( x ),~~~ x \in \Omega $ \\

with ${\cal C}^\infty$-smooth $\psi,~ \chi : \Omega \longrightarrow \mathbb{R}^a$ and $\omega
\longrightarrow \mathbb{R}$, while ${p_{\, \beta,\, \rho}} \in \mathbb{N}^n,~ l_{\, \beta,\,
\rho} \in \mathbb{N}$, and $T_{\, \beta,\, \rho} ( x, D )$ are arbitrary ${\cal
C}^\infty$-smooth coefficient polynomial nonlinear partial differential operators in $\psi$
and $\chi$. \\

{\bf Proposition 7.1.} \\

Every (MH) type nonlinear system of PDEs is resoluble.

\hfill $\Box$ \\

The main result, which extends Theorem 6.1., is the following \\

{\bf Theorem 7.1.} \\

The function $U$ in (4.2) is a global solution on $\Omega$ of the resoluble nonlinear system
of PDEs (2.1), if an only if the following {\it junction conditions} are satisfied for each
$1 \leq \beta \leq b$, in the neighbourhood of the singularity set $\Gamma$, namely \\

(7.3)~~~ $ \sum_{\,1 \,\leq\, \rho \,\leq\, r_\beta}~ \left [\, T_{\, \beta,\, \rho} ( x, D )~
                          \left (\, U_- \,, U_+ - U_- \,\right ) \,\right ]
                             D^{p_{\, \beta,\, \rho}} H_\gamma ~=~ f_\beta $

\hfill $\Box$ \\

The junction conditions (7.3) in Theorem 7.1. contain as a particular case those in Theorem
6.1., and obviously have a more simple, compact form. In the next section, indications are
given on the way these junction condition (7.3) can effectively be computed. \\ \\

{\bf 8. Computation of the Junction Conditions} \\

The relations which are useful for the effective computation of the junction conditions in
Theorem 7.1. are the following \\

(8.1)~~~ $ D^p H_\gamma ~=~ ( D^p \gamma )\, \delta_\gamma,~~~
                                          p \in \mathbb{N}^n,~~ |\, p \,| ~=~ 1 $ \\

(8.2)~~~ $ D^p ( D^l \delta_\gamma ) ~=~ ( D^p \gamma ) ( D^{l + 1} \delta_\gamma ),~~~
                               p \in \mathbb{N}^n,~~ |\, p \,| ~=~ 1,~~ l \in \mathbb{N} $ \\

(8.3)~~~ $ \gamma \delta_\gamma ~=~ 0 $ \\

(8.4)~~~ $ \gamma D^{l + 1} \delta_\gamma +
                     ( l + 1 )  D^l \delta_\gamma ~=~ 0,~~~ l \in \mathbb{N} $ \\

where $\delta_\gamma$ is the Dirac distribution concentrated on the hyper-surface $\Gamma$. In
this regard, it is well known that (8.1), (8.2) are valid when \\

(8.5)~~~ $ grad\, \gamma ( x ) ~\neq~ 0,~~~ x \in \Gamma $ \\

In view of (8.1), (8.2), one easily obtains \\

(8.6)~~~ $ D^p H_\gamma ~=~ \sum_{\, 0 \,\leq\, l \,\leq\, |\, p \,| - 1}
                              \left [~ K_{\, p,\, l} ( x, D ) \gamma ~\right ]
                             D^l \delta_\gamma,~~~ p \in \mathbb{N}^n,~ |\, p \,| \geq 1 $ \\

where $K_{\, p,\, l} ( x, D )$ are polynomial nonlinear partial differential operators with
${\cal C}^\infty$-smooth coefficients and order at most $l + 1$, which can be obtained from
the recurrent relations \\

(8.7)~~~ $ K_{\, p,\, 0} ( x, D ) \gamma ~=~ D^p \gamma,~~~
                                      p \in \mathbb{N}^n,~ |\, p \,| = 1 $ \\

(8.8)~~~ $ \begin{array}{l}
                K_{\, p + q,\, l} ( x, D ) \gamma ~=~
                  D^q \left[~ K_{\, p,\, l} ( x, D ) \gamma ~\right ]\, +
                    \,\left [~ K_{\, p,\, l - 1} ( x, D ) \gamma ~\right ]
                       D^q \gamma, \\ \\
                  \hfill ~~~~~~~~~~~~~~~~~~~~~~~~~~~~~ p, q \in \mathbb{N}^n,~ |\, q \,| = 1
            \end{array} $ \\

Now in view of (8.6), the junction conditions (7.3) can further be simplified. \\ \\

{\bf 9. Reduced Smoothness} \\

As shown in Rosinger [6, pp. 121-162], the above results hold for the more general case of
polynomial nonlinear systems of PDEs (2.1) when the coefficients are merely {\it continuous}.
In such a case, however, one has to make use not of a suitable single differential algebra of
generalized functions, but instead, of {\it chains} of such algebras. \\ \\

{\bf 10. Too Many ODEs and PDEs ?} \\

In Rubel it was shown that there exist nontrivial {\it fourth order} polynomial nonlinear
ODEs \\

(10.1)~~~ $ P ( U ( x ), U\,' ( x ), U\,'' ( x ), U\,''' ( x ), U\,'''' ( x ) ) ~=~ 0,
                                           ~~~ x \in \mathbb{R} $ \\

where $P$ is an integer coefficient polynomial in four variables, with the following {\it
dense solution} property. \\

Given any two continuous functions $f,~ w \in {\cal C}^0 ( \mathbb{R} )$, with $w ( x ) > 0$,
for $x \in \mathbb{R}$, then there exists a solution $U \in {\cal C}^\infty ( \mathbb{R} )$ of
(10.1), such that \\

(10.2)~~~ $ |~ f ( x ) \,-\, U ( x ) ~| ~\leq~ w ( x ),~~~ x \in \mathbb{R} $ \\

Two further related facts can be mentioned. The ODEs (10.1) can be constructed in a rather
simple manner. However, they are implicit in the highest derivative $U\,''''$. \\

A corresponding version of the above {\it dense solution} property was obtained for PDEs in
Buck. Namely, for every $n \in \mathbb{N},~ n \geq 2$, there exist nontrivial polynomial
nonlinear PDEs \\

(10.3)~~~ $ \begin{array}{l}
              P ( U ( x ),~\ldots~, D^p U ( x ),~\ldots ) ~=~ 0, \\ \\

              \hfill ~~~~~~~~~~~~~~~~~~~~p \in \mathbb{N}^n,~ |\, p \,| \leq m ( n ),~
                                           x \in [\, 0, 1 \,]^n
             \end{array} $ \\

with $P$ a real coefficient polynomial of degree at most $d ( n )$, such that for every
continuous function $f \in {\cal C}^0 (\, [\, 0, 1 \,]^n \,)$ and $\epsilon > 0$, there exists
a solution $U \in {\cal C}^\infty (\, [\, 0, 1 \,]^n \,)$ of (10.3), with the property \\

(10.4)~~~ $ |~ f ( x ) \,-\, U ( x ) ~| ~\leq~ \epsilon,~~~ x \in [\, 0, 1 \,]^n $ \\

This time, however, the proof is rather involved, as it uses Kolmogorov's solution to
Hilbert's Thirteenth Problem. Also, the degrees $d ( n )$ of the polynomials $P$, and the
orders $m ( n )$ of the PDEs (10.3) are very large, even for small values of $n$, such as $n =
2,\, 3,\, 4$. \\ \\

\end{document}